\documentclass[12pt]{article}
\usepackage{theorem}
\usepackage{mathrsfs}
\usepackage{amsmath}
\usepackage{amssymb}
\usepackage{epsfig}
\usepackage{graphics}

\newtheorem{thm}{Theorem}[section]
\newtheorem{prop}[thm]{Proposition}

{\theorembodyfont{\rm}
\newtheorem{defn}[thm]{Definition}
\newtheorem{prop-def}[thm]{Definition}

\newtheorem{lem}[thm]{Lemma}
\newtheorem{rem}[thm]{Remark}

}

\newcommand{\ra}{\rightarrow}

\newcommand{\fS}{\mathfrak{S}}
\newcommand{\Spec}{\mathrm{Spec\ }}
\newcommand{\Proj}{\mathrm{Proj\ }}

\newcommand{\Pic}{\mathrm{Pic}}
\newcommand{\supp}{\mathrm{supp}}
\newcommand{\GL}{\mathrm{GL}}
\newcommand{\SL}{\mathrm{SL}}
\newcommand{\Aut}{\mathrm{Aut}}
\newcommand{\bP}{{\mathbb P}}
\newcommand{\bA}{{\mathbb A}}

\newcommand{\bG}{{\mathbb G}}

\newcommand{\bQ}{{\mathbb Q}}

\newcommand{\bZ}{{\mathbb Z}}
\newcommand{\cC}{{\mathcal C}}

\newcommand{\cE}{{\mathcal E}}

\newcommand{\cL}{{\mathcal L}}
\newcommand{\cM}{{\mathcal M}}
\newcommand{\cP}{{\mathcal P}}
\newcommand{\cO}{{\mathcal O}}

\newcommand{\cX}{{\mathcal X}}
\newcommand{\oA}{\overline{A}}
\newcommand{\oM}{\overline{M}}
\newcommand{\ocM}{\overline{\mathcal M}}
\newcommand{\rM}{{\overline{\mathcal M} }^{<\iota >}_2}
\newcommand{\iHom}{\mathscr{H}\!\mathit{om}}
\newcommand{\iExt}{\mathscr{E}\!\mathit{xt}}

\date{August 2004}
\title{Classical and minimal models of the moduli space
of curves of genus two}
\author{Brendan Hassett \thanks{
	Partial support was provided by 
	National Science Foundation grants 0196187 and 0134259
	and the Sloan Foundation.}}
\begin{document}
\maketitle

\section{Introduction}

This paper is an introduction to the 
minimal model program, as applied to the moduli space of curves.
Our long-term goal is a geometric description of the canonical model
of the moduli space when it is of general type.  This entails
proving that the canonical model exists and interpreting it as
a parameter space in its own right.  

Work of Eisenbud, Harris,
and Mumford shows that $M_g$ is of general type when $g\ge 24$ 
(see \cite{HaMu} and subsequent papers).  A standard conjecture of
birational geometry--the finite generation of the canonical ring--would
imply the canonical model is 
$$\Proj \oplus_{n\ge 0} \Gamma(\oM_g,nK_{\oM_g}).$$ 
Unfortunately, this has yet to be verified in a single genus!

There is some cause for optimism:  Shepherd-Barron \cite{SB} has recently
shown that the canonical model of the moduli space of principally
polarized abelian surfaces of dimension $g\ge 12$ is 
the first Voronoi compactification.  

Another possible line of attack is to consider {\em log} canonical models of
the moduli space.  The moduli space
is best regarded as a pair $(\oM_g,\Delta)$, where $\oM_g$
is Deligne-Mumford compactification by stable curves and $\Delta$
is its boundary.  It is implicit in the work of Mumford \cite{Mu1}
that the moduli space
of stable curves is its own log canonical model (see Theorem~\ref{thm:Mumford}).
Our basic strategy it to interpolate between the log canonical model
and the (conjectural) canonical model by considering 
$$\Proj \oplus_{n\ge 0} \Gamma(\oM_g,n(K_{\oM_g}+\alpha\Delta)),$$
where $\alpha \in \bQ\cap [0,1]$ is
chosen so that $K+\alpha\Delta$ is effective.  

This program is a subject of ongoing work, inspired by correspondence
with S. Keel, in collaboration with D. Hyeon.  Future papers will address
the stable behavior of these
spaces for successively smaller values of $\alpha$.  It is remarkable
that their behavior is largely independent of the genus
(see, for example, Remark~\ref{rem:fullforce}.)  

However, for small values of $g$ special complexities arise.  
When $g=2$ or $3$, the locus in $\oM_g$ of
curves with automorphism has codimension $\le 1$.  
To include these spaces under our general framework,
we must take into account the properties of the moduli
stack $\ocM_g$.  In particular, it is necessary to 
use the canonical divisor of the moduli stack rather than
its coarse moduli space.  These differ substantially, as
the natural morphism $\ocM_g \ra \oM_g$ is ramified at stable
curves admitting automorphisms.  

Luckily, we have inherited a tremendously rich
literature on curves of small genus.  The invariant-theoretic
properties of $M_2$ were extensively studied by the 19th century
German school \cite{Cl}, who realized it as an open subset of 
the weighted projective space $\bP(1,2,3,5)$.
Theorem~\ref{thm:main} reinterprets this classical construction
using the modern language of stacks and minimal models.

\

We work over an algebraically closed field $k$ of characteristic zero. 
We use the notation $\equiv$ for $\bQ$-linear equivalence
of divisors.  
Throughout, a {\em curve} is a connected, projective, reduced scheme
of dimension one.  The {\em genus} of a curve is its arithmetic genus.  

The moduli stack of smooth (resp. stable) curves of genus $g$
is denoted $\cM_g$ (resp. $\ocM_g$);  the corresponding
coarse moduli scheme is denoted $M_g$ (resp. $\oM_g$).  
The boundary divisors in $\oM_g$ (resp. $\ocM_g$) are denoted
$$\Delta_0,\Delta_1,\ldots,\Delta_{\lfloor g/2 \rfloor} 
(\text{resp. }
\delta_0,\delta_1,\ldots,\delta_{\lfloor g/2 \rfloor}).$$
Let $Q:\ocM_g \ra \oM_g$ denote the natural morphism from the moduli
stack to the coarse moduli space,
so that
$$Q^*\Delta_i=\delta_i, i\neq 1 \quad \quad
Q^*\Delta_1=2\delta_1.$$
We write $\delta=\sum_{0\le i \le g/2}\delta_i$;  abusing
notation, we also use $\delta$ for the corresponding divisor
$$\Delta_0+1/2\Delta_1+\Delta_2+\ldots +\Delta_{\lceil g/2\rceil}$$
on $\oM_g$.

\

\noindent {\bf Acknowledgments:}  The key ideas
underlying this program were worked out in correspondence
with S. Keel.  They have been further
developed in collaboration with 
D. Hyeon.  The author also benefited from conversations
with James Spencer about rigidification, quotient stacks,
and moduli of curves of genus two, and from comments
on the manuscript by Michael van Opstall.  
Part of this paper was prepared during a visit to the
Mathematisches Institut of the Georg-August-Universit\"at,
G\"ottingen.   

\section{Classical geometry}
\subsection{Elementary facts about curves of genus two}
\label{subsect:element}

We recall results from standard textbooks, e.g.,
\cite{Ha} IV Ex 2.2  and \S 5.
Let $C$ denote a smooth curve of genus two with sheaf of differentials 
$\omega_C$.
The global sections of $\omega_C$ give the canonical
morphism
$$j:C \ra \bP(\Gamma(C,\omega_C))\simeq \bP^1$$
which is finite of degree two.  The corresponding covering
transformation $\iota:C\ra C$ is called the hyperelliptic
involution.  By the Hurwitz formula, 
$j$ is branched over six distinct points
$$\{b_1,\ldots,b_6 \} \subset \bP(\Gamma(C,\omega_C)).$$
On fixing an identification 
$\bP(\Gamma(C,\omega_C))\simeq \bP^1$, we can write down
a nontrivial binary sextic form vanishing at the branch points
$$
F \in \Gamma(\bP^1,\cO_{\bP^1}(6)),$$
determined by 
$\{b_1,\ldots,b_6\} \in \bP^1$ up to a scalar.  

Conversely, suppose we have a binary sextic form $F$
with six distinct zeros $b_1,\ldots, b_6 \in \bP^1$.  Then there is a 
unique degree-two cover of $\bP^1$ branched over these points.  This 
is a smooth curve of genus two and the map to $\bP^1$
is the canonical morphism.  Moreover, the isomorphism
class of $C$ depends only on the orbit of $F$ under the action of 
$\GL_2$.  

To summarize: There is a one-to-one correspondence between 
isomorphism classes of curves of genus two and $\GL_2$-orbits
of binary sextic forms with distinct zeros.  

\

We will need a relative version of this dictionary,
following \cite{Vi}.  Let $\pi:\cC \ra S$ be a smooth
morphism to a scheme of finite type over $k$, with fibers
curves of genus two.  Since the relative dualizing sheaf $\omega_{\pi}$
is globally generated, there is a relative double cover
$$
\begin{array}{rcccl}
\cC & & \stackrel{j}{\longrightarrow} & &  
		\bP^1(\pi_*\omega_{\pi}):=\bP \\
   &{\scriptstyle \pi} \searrow & & \swarrow {\scriptstyle \psi} &  \\
   & & S & & 
\end{array}
$$
with associated involution $\iota$.  
Using the trace we decompose 
$j_*\cO_{\cC}\simeq \cO_{\bP}\oplus \cL$, where 
$\cL$ has relative
degree $-3$ on the fibers of $\psi$.  The $\cO_{\bP}$-algebra structure
on $\cO_{\cC}$ is thus determined by an isomorphism 
$\cL^2\ra \cO_{\bP}$, i.e., by a nonvanishing section
of $\cL^{-2}$ with zeros along the branch locus of $j$.
By relative duality
$$j_*\omega_{\pi}=\iHom_{\cO_{\bP}}(j_*\cO_{\cC},\omega_{\psi})
\simeq \omega_{\psi} \oplus (\omega_{\psi}\otimes \cL^{-1})$$
which yields
$$\pi_*\omega_{\pi} \simeq \psi_*(\omega_{\psi}\otimes \cL^{-1}).$$
Using the identifications 
$$\psi_*\cO_{\bP}(+1)=\pi_*\omega_{\pi}, \quad 
\omega_{\psi}=(\psi^*\det \pi_*\omega_{\pi})(-2), \quad
\Pic(\bP)=\Pic(S)\oplus \bZ c_1(\cO_{\bP}(+1)),$$
we find
$$\cL^{-1}\simeq \cO_{\bP}(+1)\otimes \omega_{\psi}^{-1}=
\cO_{\bP}(3)\otimes (\psi^* \det \pi_*\omega_{\pi})^{-1}.$$
The class of the branch divisor is thus
\begin{equation} \label{eq:Chern}
-2c_1(\cL)=c_1(\cO_{\bP}(6))-2\psi^* c_1( \det \pi_*\omega_{\pi}).
\end{equation}

This has practical implications:
An isomorphism of $C$ induces a linear transformation on $\Gamma(C,\omega_C)$,
which respects the binary sextic $F$ up to a scalar.  
Formula (\ref{eq:Chern}) allows us to keep track of this scalar.
For each $M\in \GL_2$, we have the linear action
$$(x,y) \mapsto (x,y)\left( \begin{matrix} m_{11} & m_{12} \\
					   m_{21} & m_{22}
			\end{matrix} \right)$$
which induces a natural left action on binary sextic forms
$$F \mapsto 
F(xm_{11}+ym_{21}, xm_{12}+ym_{22}).$$
We normalize this action using
formula (\ref{eq:Chern}) 
\begin{equation} \label{eq:action}
F \mapsto (M,F):=(\det M)^{-2}F(xm_{11}+ym_{21}, xm_{12}+ym_{22}),
\end{equation}
so that $(M,F)=F$ if and only $M$ is 
induced from an automorphism of $C$.  

\

A smooth curve 
is {\em bielliptic} if it admits a degree-two morphism
$i: C\ra E$ to an elliptic curve;
the covering transformation is called a {\em bielliptic involution}.
For curves of genus two any bielliptic involution
commutes with the hyperelliptic involution, which yields a diagram
$$
\begin{array}{ccc}
C & \stackrel{i} \ra & E \\
{\scriptstyle j} \downarrow \quad & & 
		\quad \downarrow  
{\scriptstyle \bar{j}} \\
\bP^1 & \stackrel{\bar{i}}\ra & \bP^1
\end{array}
$$
where $\bar{i}$ and $\bar{j}$ are the double
covers induced on
quotients.  The branch locus of $j$ is preserved
by the covering transformation for
$\bar{i}$, which is conjugate to 
$[x,y] \mapsto [y,x]$.  
The resulting involution of the
branch locus will also be
called a bielliptic involution.
Thus $C$ is isomorphic to a double cover branched over
$$\{[\alpha_1,1],[1,\alpha_1],[\alpha_2,1,][1,\alpha_2],[\alpha_3,1],
[1,\alpha_3]\}$$
for some $\alpha_1,\alpha_2,\alpha_3\in k$.  Conversely, each such curve
admits a diagram as above and thus is bielliptic.

\subsection{Invariant theory of binary sextics}
\label{subsect:invariant}
We observe the classical convention for normalizing the
coefficients of a binary sextic
$$
F=ax^6+6bx^5y+15cx^4y^2+20dx^3y^3+15ex^2y^4+6fxy^5+gy^6.$$ 
The action (\ref{eq:action}) induces an action of 
$\GL_2$ on $k[a,b,c,d,e,f,g]$.
Recall that 
a polynomial $P \in k[a,b,c,d,e,f,g]$ is 
{\em $\SL_2$-invariant} if 
for each $M\in \SL_2$,
we have
$$(M,P) = P.$$
We write 
$$R:=k[a,b,c,d,e,f,g]^{\SL_2}$$ 
for the ring of such
invariants.  
If $P$ is $\SL_2$-invariant then each homogeneous component of $P$
is as well, so $R$ is a graded ring.  Every
homogeneous invariant satisfies the functional
relation
\begin{equation} \label{eq:functional}
(M,P)=(\det M)^{\deg(P)}P, \quad M \in \GL_2;
\end{equation}
here it is essential that the action 
(\ref{eq:action}) include the factor $(\det M)^2$.  
The transformation
$$(x,y) \ra (y,x)$$
thus reverses the sign of invariants of odd degree.  These
are called {\em skew invariants} in the classical literature.  

Explicit generators for $R$ were first written down in the
nineteenth century, e.g., \cite{Cl}, pp. 296, in symbolic notation,
\cite{Ca} and \cite{Sa} as explicit polynomials--the 
second edition of Salmon's {\em Higher algebra} has the
most detailed information, and also
\cite{El} pp. 322.
A nice early twentieth-century discussion is \cite{Sc} pp. 90
and a modern account invoking the representation theory
of $\SL_2$ is \cite{Sp}.  

For our purposes,
the symmetric function representation of the invariants in
\cite{Ig2} pp. 176 and 185 is the most useful.  Let $\xi_1,\ldots,\xi_6$
denote the roots of the dehomogenized form $F(x,1)$, and
write $(ij)$ as shorthand for $\xi_i-\xi_j$.  We write
\begin{eqnarray*}
A&=&a^2 \sum_{\text{fifteen}} (12)^2(34)^2(56)^2 \\
B&=&a^4 \sum_{\text{ten}} (12)^2(23)^2(31)^2(45)^2(56)^2(64)^2 \\
C&=&a^6 \sum_{\text{sixty}} (12)^2(23)^2(31)^2(45)^2(56)^2(64)^2
		(14)^2(25)^2(36)^2 \\
D&=& a^{10}\prod_{ij}(ij)^2 \\
E&=& a^{15}\prod_{\text{fifteen}} \det \left(
\begin{matrix} 1 & \xi_1+\xi_2 & \xi_1\xi_2 \\
 1 & \xi_3+\xi_4 & \xi_3\xi_4 \\
 1 & \xi_5+\xi_6 & \xi_5\xi_6 
\end{matrix}
\right) = a^{15} \prod_{\text{fifteen}} \left(
(14)(36)(52)-(16)(32)(54)\right)
\end{eqnarray*}
where the summations are chosen to make the 
expressions $\mathfrak{S}_6$-symmetric.  Consequently, 
$A,B,C,D,$ and $E$ can all be expressed 
as polynomials in $\bQ[a,b,c,d,e,f,g]$, e.g.,
\begin{eqnarray*}
A&=&-240(ag-6bf+15ce-10d^2)\\
B&=&-162000\det \left(\begin{matrix} a & b & c & d \\
				      b & c & d & e \\
				      c& d & e & f \\
				      d & e & f & g 
			\end{matrix} \right) + 1620 (ag-6bf+15ce-10d^2)^2.
\end{eqnarray*}
In classical terminology, $(ag-6bf+15ce-10d^2)$ is the {\em sixth
transvectant of $F$ over itself};  transvection is one of the main 
operations in Gordan's proof of finiteness for invariants of binary forms.
The determinantal
expression is the {\em catalecticant} of $F$:  It vanishes
precisely when $F$ can be expressed as a sum of three sixth
powers \cite{El} pp. 276.  

The following facts will
be useful for subsequent analysis:
\begin{prop} \label{prop:inv}
\begin{enumerate}
\item{The expressions $A,B,C,D,$ and $E$ are invariant 
and generate $R$ \cite{Hi}, pp. 100, \cite{Cl}, etc.}
\item{$D$ is the discriminant and vanishes
precisely when the binary form has a multiple root.}
\item{$B,C,D,$ and $E$ vanish whenever the binary form has a triple root;
$A$ vanishes when the form has a quadruple root.}
\item{$E$ vanishes if and only if the form admits a bielliptic involution, as defined
in \S\ref{subsect:element}
\cite{El}, pp. 327 and \cite{Cl}, pp. 457.}
\item{The unique irreducible relation among the invariants is
$$E^2=G(A,B,C,D),$$
where $G$ is weighted-homogeneous of degree $30$ \cite{Cl}, pp. 299.}
\end{enumerate}
\end{prop}
The notation used for the generating invariants is
not consistent among authors.  Our notation is
consistent with that of Igusa, but inconsistent with Clebsch's
and Salmon's.  Of course, the invariants of degree two and fifteen
are unique up to scalar.

\subsection{The projective invariant-theory quotient}
We consider 
$$X:=\Proj R= \Proj \frac{k[A,B,C,D,E]}{\left<E^2-G(A,B,C,D)\right>}.$$
If $A=B=C=D=0$ then $E=0$ as well, so $X$ is covered by the distinguished
affine open subsets
$$\{A\neq 0 \} \quad \{B \neq 0 \} \quad \{C \neq 0 \}
\quad  \{ D\neq 0 \}.$$
However, in each localization
$$(k[A,B,C,D,E][A^{-1}])_0 \quad
(k[A,B,C,D,E][B^{-1}])_0$$
$$(k[A,B,C,D,E][C^{-1}])_0 \quad 
(k[A,B,C,D,E][D^{-1}])_0
$$
only even powers of $E$ appear, so all the
functions over these distinguished open subsets can be expressed
in terms of $A,B,C,D$.  In light of 
Proposition~\ref{prop:inv}, we find
\begin{prop} \label{prop:weight1}
\begin{enumerate}
\item{$X\simeq
\Proj k[A,B,C,D] \simeq \bP(2,4,6,10)\simeq \bP(1,2,3,5)$
\cite{Ig2}, pp. 177.}
\item{A binary sextic with a zero of multiplicity three, admitting
a nonvanishing invariant of positive degree, is
mapped to $p:=[1,0,0,0,0] \in X$.}
\item{All positive-degree invariants vanish at binary sextics with a 
zero of multiplicity four;  they do not yield
points of $X$. } 
\end{enumerate}
\end{prop}

Geometric Invariant Theory gives an interpretation
of the points of $X$:
\begin{prop} \label{prop:GIT}
\begin{enumerate}
\item{a binary sextic is stable (resp. semistable) if and only if 
its zeros have multiplicity $\le 2$ (resp. $\le 3$)
\cite{GIT}, ch. 4 \S 1;}
\item{$X - \{[1,0,0,0,0]\}$ is a geometric quotient for binary
sextics with zeros of multiplicity $\le 2$ \cite{GIT} 1.10.}
\end{enumerate}
\end{prop}
The `only if' part of the first assertion can be deduced from 
Proposition~\ref{prop:weight1}.
As $X-\{D=0\}$ is a geometric quotient
for binary sextics with distinct zeros, 
our analysis of genus two curves in 
\S\ref{subsect:element} yields
\begin{prop}\label{prop:m2}
The moduli scheme $M_2$ can be identified
with 
$X-\{D=0\}$,
where $D$ is the discriminant.
\end{prop}
\begin{rem}
This construction 
definitely fails in characteristic two.  
If the double cover $j:C\ra \bP^1$
is wildly ramified, the branch divisor may have
multiplicities $>3$.  These curves correspond
to unstable points under the $\SL_2$-action, and
thus are not represented in the invariant-theory
quotient.  \cite{Ig1} has a detailed account of what must be
done in this case.  
\end{rem}

\subsection{Invariant-theory quotient as a contraction}
We sketch the relationship between the invariant-theory
quotient and the moduli space of stable curves.  
\begin{defn}
A birational map of normal projective varieties 
$$\beta:Y \dashrightarrow  X$$
is a {\em contraction} if $\beta^{-1}$ has no exceptional
divisors, i.e., the proper transform of each codimension-one
subset in $X$ has codimension one in $Y$.
\end{defn}

\begin{prop} \label{prop:bircont}
There exists a birational contraction $\beta:\oM_2\dashrightarrow X$
restricting to the identity along the open
subset $M_2$.  $\beta$ 
is an isomorphism over $\oM_2 - \Delta_1$ and
contracts $\Delta_1$
to the point $p$.  
\end{prop}
{\em proof:} To produce the birational contraction,
we exhibit a morphism
$$\beta^{-1}:U \hookrightarrow \oM_2$$
where $U\subset X$ is open with complement of codimension $\ge 2$
and $\beta^{-1}|_{M_2\cap U}$ is the identity.  
We shall take $U=X-p$, where $p$ corresponds
to the binary forms with a triple zero (cf. Proposition~\ref{prop:weight1}.)  

The universal binary sextic is a hypersurface
$$
W:=\{ax^6+6bx^5y+15cx^4y^2+20dx^3y^3+15ex^2y^4+6fxy^5+gy^6=0\}
\subset \bA^7 \times \bP^1.
$$
Its class in $\mathrm{Pic}(\bA^7 \times \bP^1)$
is divisible by two, so
there exists a double cover $\cC'\ra \bA^7\times \bP^1$ 
simply branched over $W$.  Composing with the projection onto
the first factor, we obtain a morphism
$$\pi':\cC' \ra \bA^7.$$
Let $S\subset \bA^7$ denote the open subset corresponding to forms
whose zeros all have multiplicity $\le 2$ and
$$\pi:\cC \ra S$$
the restriction of $\pi'$ to $S$.  Since $\pi$ is a composition
of flat morphisms, it is also flat.

Consider the fiber of $\pi$ over a given
binary sextic $F$:  It is a double cover $j:\cC_F\ra \bP^1$ 
branched over the zeros of $F$.  We claim $\cC_F$ is a stable
curve of genus two, not contained in $\Delta_1$.  
Evidently $\cC_F$ is smooth and simply branched over
the zeros with multiplicity one.
Over the double zeros
$\cC_F$ has local equation $y^2=x^2$,
which defines a node.  We have
$j^*\cO_{\bP^1}(+1)=\omega_{\cC_F}$, which is ample on $\cC_F$, 
so $\cC_F$ is stable.  
The normalization
$\nu:\cC^{\nu}_F\ra \cC_F$ is the double cover branched along $F'=0$,
where $F'$ is the product of the factors of $F$ with
multiplicity one;  $\cC_F$ is obtained from $\cC^{\nu}_F$ by gluing
the pairs of points over the each double root of $F$.
There are three possibilities:
\begin{enumerate}
\item{$\deg(F')=4$, in which case
$\cC^{\nu}_F$ is connected of genus one;}
\item{$\deg(F')=2$, in which case
$\cC^{\nu}_F$ is connected of genus zero;}
\item{$\deg(F')=0$, in which case
$\cC^{\nu}_F$ has two connected components of genus zero.}
\end{enumerate}
Since $\cC_F$ cannot be expressed as the union of two
subcurves of genus one meeting at a point, the resulting curve
is not in $\Delta_1$.  

The classifying morphism
$$S\ra \ocM_2-\delta_1$$
is equivariant with respect to the $\GL_2$-action on binary sextics,
and therefore descends to a morphism
$$\Phi:S/\GL_2 \ra \ocM_2-\delta_1.$$
We remark that this is a morphism of stacks.
Since $U$ is a geometric quotient for binary sextics with
zeros of multiplicity $\le 2$ (see Proposition~\ref{prop:GIT}),
$U$ is also the coarse moduli space for $S/\GL_2$.
We define $\beta^{-1}$ to be the induced morphism on coarse moduli spaces. 

It remains to show this is bijective onto its image.
Suppose we 
are given a stable curve $C$ of genus two not contained
in $\Delta_1$.  Quite generally, $\omega_C$ is globally generated
for any stable curve without disconnecting nodes;  the only
curves in $\oM_2$ with disconnecting nodes lie in $\Delta_1$.  
Thus the sections of $\omega_C$ give a double cover
$$j:C \ra \bP^1$$
branched along a sextic, with zeros of multiplicity $\le 2$
because $C$ is nodal.  The analysis above shows that every
such sextic arises in this way.  $\square$

\subsection{Blowing up the invariant-theory quotient}
We recall the principal result of \cite{Ig2}.  
Let $A_g$ denote the moduli space of principally polarized
abelian varieties of dimension $g$, $\oA_g$ its Satake compactification.
Recall that $\oA_g=\Proj S$, where $S$ is the ring
of $\mathrm{Sp}(g,\bZ)$-modular forms;  we use
$\lambda$ to denote the resulting polarization on $\oA_g$.  Let
$t:M_g \hookrightarrow A_g$ denote
the Torelli morphism, associating to each curve its Jacobian.  

Now assume $g=2$.  Regarding $M_2$ as an open subset of
$X$ (see Proposition~\ref{prop:m2}),  $t$ extends to a rational map
$$\tau:X \dashrightarrow \oA_2.$$
The inclusion 
and the Torelli morphism induce
\begin{equation}
\label{eq:compact}
M_2 \hookrightarrow \widetilde{X}:={\overline{\mathrm{Graph}(\tau)}} \subset 
X \times \oA_2.
\end{equation}
In particular, $\widetilde{X}$ compactifies $M_2$.  

\begin{thm}\label{thm:Igusa}\cite{Ig2}
The indeterminacy of $\tau$ is the point $p=[1,0,0,0,0]\in X$ corresponding
to binary sextic forms with a zero of multiplicity three.  If we choose
local coordinates at $p$
$$x_1=2^43^2B/A^2 \quad x_2=2^63^3(3C-AB)/A^3 \quad x_3=2\cdot3^5D/A^5$$
then $\tau$ is resolved by a weighted blow-up centered at $p$
$$b:\widetilde{X} \ra X$$
with weights
$$\text{weight}(x_1)=2 \quad
\text{weight}(x_2)=3 \quad
\text{weight}(x_3)=6.$$
The exceptional divisor of $b$ is mapped isomorphically
to the locus of principally
polarized abelian surfaces that decompose as a product of two
elliptic curves (with the induced product polarization).  
\end{thm}
\begin{rem}\label{rem:Igusa}
Igusa's result is considerably more precise:  
He explicitly computes the correspondence between the ring of
invariants $R$ and the ring of modular forms $S$.  
In particular, $S$ is a polynomial ring with
generators in degrees $4,6,10,12$ and the locus of products
is given by the vanishing of a form of weight $10$.
\end{rem}

\subsection{Comparing the blow-up with moduli space}
\label{subsect:comp}
\begin{prop}\label{prop:codimone}
The open imbedding $M_2 \hookrightarrow \widetilde{X}$ extends
to a birational map
$$\gamma: \oM_2 \dashrightarrow \widetilde{X}$$
which is an isomorphism in codimension one.
\end{prop}
In particular, $\gamma$ and
$\gamma^{-1}$ are both birational contractions.  
In Proposition~\ref{prop:nef2}  
we will prove that $\gamma$ is an isomorphism.

\noindent {\em proof:}
The Torelli morphism admits an extension
$\overline{t}: \oM_g \ra \oA_g$
\cite{Na1}, Theorem 3.
This is not an isomorphism for $g>1$:  
The divisor $\Delta_0 \subset \oM_g$ 
is mapped to a boundary stratum of $\oA_g$, which has
codimension $\ge 2$.  However, in genus two 
$\overline{t}$ is
an isomorphism at the generic point of $\Delta_1$.  
Indeed, the Jacobian of a curve
$[E_1\cup_q E_2] \in \Delta_1$,
with $E_1$ and $E_2$ smooth of genus one,
is the abelian surface $E_1\times E_2$.  

The following diagram
summarizes the various birational maps and morphisms:
$$
\begin{array}{ccc}
\oA_2   &\stackrel{\pi}{\leftarrow}  & \widetilde{X} \\
{\scriptstyle \overline{t}} \uparrow \quad  &  &
        \quad  \downarrow {\scriptstyle b} \\
\oM_2 &\stackrel{\beta}{\dashrightarrow} & X
\end{array}.
$$
By Theorem~\ref{thm:Igusa}, $\pi$ is also an isomorphism 
over the generic point
of $\overline{t}(\Delta_1)$, so $\gamma$ is an isomorphism at the generic
point of $\Delta_1$.  $\beta$ and $b$ are both isomorphisms over 
the generic point of the divisor $\beta(\Delta_0)$ 
(see Proposition~\ref{prop:bircont}), so $\gamma$ is an isomorphism
at the generic point of $\Delta_0$.  Since $\gamma$ is regular along
$M_2=\oM_2-\Delta_0-\Delta_1$, the result follows.  $\square$

The proper transforms of $\Delta_0$ and $\Delta_1$ in $\widetilde{X}$
are denoted $\widetilde{\Delta}_0$ and $\widetilde{\Delta}_1$.  
Thus $\widetilde{\Delta}_1$ is the exceptional divisor of
$b:\widetilde{X} \ra X$.  

\begin{rem}[Bibliographic note]
There are a
number of partial desingularizations of $\oA_g$ through which
$\overline{t}$ factors, e.g., 
the `Igusa monoidal transform'
\cite{Ig3}, \cite{Na1}
and the toroidal compactification associated to the 2nd Voronoi
fan \cite{Na2}.   When $g=2$, these approaches
coincide \cite{Na2} Remark 2.8 and yield a partial 
desingularization
$\widetilde{A}_2 \ra \oA_2.$
See \cite{Ig3} Theorem 5 for a blow-up representation, expressed in
terms of modular forms;  the center of this blow-up is 
in the boundary $\oA_2 - A_2$.  
Namikawa \cite{Na1} \S 9 has shown that the factorization
$\oM_2 \stackrel{\sim}{\ra} \widetilde{A}_2$
is an isomorphism.  
Notwithstanding
Igusa's explicit formulas for $\tau:X \dashrightarrow \oA_2$
and $\widetilde{A}_2 \ra \oA_2$, it is not entirely obvious 
how to extract an isomorphism
$\widetilde{A}_2 \stackrel{\sim}{\ra} \widetilde{X}.$
\end{rem}

\section{Stack geometry}
\subsection{A stack-theoretic quotient}
Proposition~\ref{prop:weight1} might suggest
that the invariant $E$ is irrelevant to the geometry of the 
quotient.  However, we have so far ignored possible
{\em stack structures} on the quotient, which are intertwined
with the geometry of $E$.  
There are a number of natural stacks to consider,
including the $\GL_2$-quotient stack.  
Our choice is dictated by pedagogical
imperatives, i.e., to exhibit concretely the nontrivial
inertia along the bielliptic locus where $E$ vanishes. 

The ring of invariants $R$ is graded by degree, so we have a
natural $\bG_m$-action on the affine variety
$Y=\Spec R$.  Now $\bG_m$ acts 
on the open subset $Y-(0,0,0,0,0)$ with finite stabilizers
and closed orbits, so the quotient stack
$$\cX:=\left( Y-(0,0,0,0,0) \right)/\bG_m$$
is a separated Deligne-Mumford stack with coarse moduli space 
$q:\cX \ra X$ (see \cite{LM}
10.13.2,7.6,8.1 for more information).  
The points of $Y-(0,0,0,0,0)$ with nontrivial stabilizer map to
the points of $\cX$ with nontrivial inertia groups;  this is 
the ramification locus of $q$. 

We collect some geometric properties of $\cX$:
\begin{prop} \label{prop:ramq}
\begin{enumerate}
\item{The closed imbedding
$\Spec R \hookrightarrow \Spec k[A,B,C,D,E] $
induces a closed imbedding
$i:\cX \hookrightarrow \cP(2,4,6,10,15)$,
where 
$$\cP(2,4,6,10,15):= \left(\Spec k[A,B,C,D,E]-(0,0,0,0,0)\right)/\bG_m,$$
with $\bG_m$ acting with weights $(2,4,6,10,15)$
$$t\cdot (A,B,C,D,E) \mapsto (t^2A,t^4B,t^6C,t^{10}D,t^{15}E);$$}
\item{the dualizing sheaf of $\cX$ is given via adjunction
$$\omega_{\cX}=i^*\omega_{\cP(2,4,6,10,15)}(\cX)\simeq 
i^*\cO_{\cP(2,4,6,10,15)}(-7),$$
where $\cO_{\cP(2,4,6,10,15)}(+1)$ is the invertible sheaf 
associated with the principal $\bG_m$-bundle
classified by the identity character of $\bG_m$;}
\item{the ramification divisor of $q$ is 
$$\cE:=\{ E=0 \} \subset \cX$$
and we have 
$$q^*\omega_X\simeq \omega_{\cX}(-\cE)\simeq i^*\cO_{\cP(2,4,6,10,15)}(-22).$$}
\end{enumerate}
\end{prop}
{\em proof:} The first two assertions do not require proof.
As for the third, it follows from the classification of 
possible automorphisms of binary sextics \cite{Bo}, \cite{Ig1} \S 8
The only automorphism type occuring in codimension one
is the bielliptic involution (cf. Proposition~\ref{prop:inv}).
$\square$

\begin{rem} Notwithstanding Proposition~\ref{prop:m2}, 
$\cM_2$ is {\em not} contained in $\cX$ as an open substack.  
Using the functional relation (\ref{eq:functional}),
the inertia group of $[F]\in \cX-p$ is the quotient
$$\{M\in \GL_2: (M,F)=F \}/\{M\in \SL_2: (M,F)=F \},$$
which is trivial for generic binary forms.  

The inertia group at $[C]\in \cM_2$ is $\Aut(C)$, and
the presence of the hyperelliptic involution $\iota$ means 
this is always nontrivial.
Now $\Aut(C)$ has a natural representation on
$\Gamma(C,\omega_C)$ and an induced representation
on $\wedge^2\Gamma(C,\omega_C)$ that is not faithful:
We do not see elements of $\Aut(C)$ acting on $\Gamma(C,\omega_C)$
with determinant one, e.g., $\iota$, which acts on 
$\Gamma(C,\omega_C)$ by $-I$.  The corresponding quotient of $\Aut(C)$
is the inertia picked up by $\cX$.  
\end{rem}

\subsection{Analysis of boundary divisors in the moduli stack}
\label{subsect:stackgeo}
\begin{lem} \label{lem:involution}
Every stable curve of genus two admits a canonical 
hyperelliptic involution, which is central in its automorphism
group.
\end{lem}
{\em proof:}
We claim that every stable curve of genus two is canonically 
a double cover of a nodal curve
of genus zero 
$$j:C \ra R, \quad 
R=\begin{cases} \bP^1 & \text{ if } [C]\not \in \Delta_1\\
            \bP^1 \cup_r \bP^1 & \text{ if } [C] \in \Delta_1
\end{cases}	
$$
branched over six smooth points $b_1,\ldots,b_6 \in R$
and the node $r\in R$.  
The covering transformation $\iota$ therefore commutes with
each automorphism of $C$.

The cover is induced by 
$$C \ra \bP(\Gamma(C,\omega^2_C))\simeq \bP^2$$
which factors
$$C \stackrel{j}{\ra} R \subset \bP^2,$$
where $R$ is a plane conic
and $j$ is finite of degree two.  Indeed, for curves not 
in $\Delta_1$ this 
is the double cover $C\ra \bP^1$ discussed in the proof of 
Proposition~\ref{prop:bircont}.  For curves $C=E_1\cup_q E_2 \in \Delta_1$
with $q=E_1\cap E_2$ the disconnecting node joining the genus-one
components $E_1$ and $E_2$, we have a 
double cover
$$j:E_1\cup_q E_2 \longrightarrow \bP^1 \cup_r \bP^1,$$
with $j(q)=r$ and $j$ mapping the genus one components two-to-one onto
rational components, with ramification at $q$ along each component.  
$\square$

Consider stable curves of genus two with automorphisms beyond 
$\iota$.  There are two possibilities:
Either $(R,b_1,\ldots,b_6)$ admits automorphisms permuting
the $b_i$ or $j:C \ra R$ admits covering transformations
other than the canonical hyperelliptic involution.  
The classification
of automorphism groups (\cite{Bo} or \cite{Ig1}\S 8) yields the following
possibilities in codimension one:
\begin{enumerate}
\item{the curves $C$ in $\Delta_1$;  here
$j:C \ra \bP^1 \cup \bP^1$
admits involutions fixing each component of $C$;}
\item{the closure of the locus of curves $j:C \ra \bP^1$
branched over six points admitting a bielliptic involution.}
\end{enumerate}

\

At each point $[C]\in \oM_2$, the moduli space
is \'etale-locally isomorphic to the quotient
$T_{[C]}\ocM_2/\Aut(C)$ at the origin, where
$$T_{[C]}\ocM_2=\mathrm{Ext}^1(\Omega^1_C,\cO_C))$$
is the tangent space with the induced automorphism action.
When $C$ is smooth, Serre duality gives
$$T_{[C]}\ocM_2=\Gamma(C,\omega^2_C)^*;$$
equation (\ref{eq:action}) shows that
the hyperelliptic
involution acts trivially and a
bielliptic involution acts by reflection. 

This local isomorphism can be chosen
so the divisors $\Delta_0$ and $\Delta_1,$ correspond
to the images of unions of distinguished hyperplanes in $T_{[C]}\ocM_2$.
The local-global spectral sequence gives an exact sequence
\begin{equation} \label{ex:locglob}
0 \ra H^1(\iHom(\Omega^1_C,\cO_C)) \ra \mathrm{Ext}^1(\Omega^1_C,\cO_C) \ra
\Gamma(\iExt^1(\Omega^1_C,\cO_C))\ra 0.
\end{equation}
A local computation implies
$$\iExt^1(\Omega^1_C,\cO_C) \simeq
\oplus_{\text{nodes }p\in C} k.$$
If $C=E_1\cup_q E_2$ is in $\Delta_1$, let
$\widehat{\delta_1}\subset T_{[C]}\ocM_2$ 
denote the hyperplane corresponding to 
the node disconnecting
the genus-one components of $C$, i.e.,
the kernel of the projection onto the direct summand
corresponding to $q$.  
The extra covering involutions of $j:C\ra \bP^1\cup \bP^1$ 
mentioned above
act on $T_{[C]}\ocM_2$,
trivially on $\widehat{\delta_1}$
and by multiplication by $(-1)$ on $\iExt^1(\Omega^1_C,\cO_C)_q$.
This corresponds to reflection across $\widehat{\delta_1}$.  
Let $\widehat{\delta_0}$ denote the union of the hyperplanes
corresponding to each of the non-disconnecting nodes of $C$.
The surjectivity of the last arrow in (\ref{ex:locglob})
means that 
$\widehat{\delta}=\widehat{\delta_1}\cup \widehat{\delta_0}$ 
is normal crossings.

\begin{defn}
Let $\xi \subset \ocM_2$ (resp. $\Xi \subset \oM_2$)
denote the closure of the 
smooth curves admitting
a bielliptic involution. 
\end{defn}
These are irreducible of codimension one, e.g., by Proposition~\ref{prop:inv}
and the characterization of bielliptic curves.  This can also be seen
infinitesimally:  Under the local identification with $T_{[C]}\ocM_2/\Aut(C)$,
each branch of $\xi$ is identified with the hyperplane 
$T_{[C]}\ocM_2$ fixed by the corresponding bielliptic
involution (acting by reflection
on the tangent space).  The union of these hyperplanes
is denoted $\widehat{\xi}$.

The bielliptic divisor 
has more complicated local
geometry, as $\Xi$ may have quite a few local
branches.  
For example, 
$$F(x,y)=xy(x-\xi y)(x+\xi y)(x - \xi^{-1}y)(x+\xi^{-1}y)$$
has automorphism group isomorphic to the Klein four-group and
admits two involutions
$$[x,y]\mapsto [y,x] \quad 
[x,y] \mapsto [-y, x].$$
The cyclotomic form
$$F(x,y)=x^6+y^6=(x-\zeta y)(x-\zeta^3 y)(x-\zeta^5 y)(x-\zeta^7 y)
(x-\zeta^9 y)(x-\zeta^{11}y) \quad \zeta^{12}=1$$
has automorphism group isomorphic to the dihedral group with $12$
elements and
admits {\em four} distinct involutions 
$$[x,y]\mapsto [y,x] \quad 
[x,y]\mapsto [-x,y] \quad 
[x,y] \mapsto [\zeta^4 y, x] \quad
[x,y] \mapsto [\zeta^8 y, x] .
$$
In particular, 
$\widehat{\xi}$
is {\em not} normal crossings.

\begin{prop}[Ramification formula for $Q:\ocM_2 \ra \oM_2$] 
\label{prop:ramify}
$$
K_{\ocM_2}+\alpha \delta \equiv Q^*(K_{\oM_2}+ \alpha \Delta_0+
\frac{1+\alpha}{2}\Delta_1+\frac{1}{2}\Xi)
$$
\end{prop}
\noindent {\em proof:}
Lemma~\ref{lem:involution} allows us to consider 
the rigidification (\cite{ACV} \S 5) of $\ocM_2$
with respect to the group $\left<\iota\right>$
generated by the canonical involution
$$r:\ocM_2 \ra \rM.$$
Given a scheme $T$, $T$-valued points of $\rM$ 
correspond to families of genus two 
stable curves over $T$, where two families are identified 
if they differ by a $\left<\iota\right>$-valued cocycle over $T$.
The inertia group of 
$\rM$ at $[C]$ is the quotient
$$\Aut(C)/\left< \iota \right>.$$
We have a factorization
$$\ocM_2 \stackrel{r}{\ra} \rM \stackrel{\bar{q}}{\ra} \oM_2$$
so that $r$ is \'etale of degree two and $\oM_2$ is the
coarse moduli space of $\rM$ \cite{ACV} Theorem 5.1.5.  
We therefore obtain the following formula for dualizing sheaves
$$r^*\omega_{\rM}\simeq \omega_{\ocM_2}.$$

Let $\delta_1^{<\iota>}$ and $\xi^{<\iota>}$ denote the
corresponding Cartier divisors in $\rM$.
As $\bar{q}$ has simple ramification along these divisors,
we obtain
\begin{eqnarray*}
\bar{q}^*\Delta_1=2\delta_1^{<\iota >}    & &
\bar{q}^*\Xi=2\xi^{<\iota >}    \\
{\bar{q}}^*K_{\oM_2}&=&K_{\rM}-\delta_1^{<\iota>}-\xi^{<\iota>},
\end{eqnarray*}
which together imply the formula.$\square$

\subsection{Analysis 
of $b:\widetilde{X}\ra X$ along the exceptional divisor}
By Theorem~\ref{thm:Igusa}, the exceptional divisor $\widetilde{\Delta}_1$
is mapped isomorphically to the locus in $\oA_2$ parametrizing
abelian surfaces
decomposing into products of elliptic curves (as a principally
polarized abelian variety);  this 
is isomorphic to $\bP(2,3,6)$.  
\begin{prop}\label{prop:symsquare}
Let $\widetilde{\Xi}$ denote the proper transform of $\Xi$
in $\widetilde{X}$.  We have the following formulas:
\begin{eqnarray*}
K_{\widetilde{X}}&\equiv& b^*K_X+10\widetilde{\Delta}_1\\
b^*\{G=0\}&\equiv& \widetilde{\Xi}
+12\widetilde{\Delta}_1\\
b^*\{D=0\}& \equiv & \widetilde{\Delta}_0
+6\widetilde{\Delta}_1 \\
\widetilde{\Xi}&\equiv&
3\widetilde{\Delta}_0+12\widetilde{\Delta}_1.
\end{eqnarray*}
\end{prop}
We pause to explore the geometry of 
$\widetilde{\Delta}_1$.  Naively, one might expect
this to be the symmetric square of the moduli space of elliptic
curves.  However, in taking symmetric squares we should be mindful of
the stack structures.
The coarse moduli space of the symmetric square need not be
isomorphic to the symmetric square of the coarse
moduli space.

The standard theory of modular forms implies 
$$\ocM_{1,1}\simeq \cP(4,6)=
\left(\Spec k[g_2,g_3] - (0,0)\right)/\bG_m,$$
with $\bG_m$ acting with weights $(4,6)$
$$t\cdot (g_2,g_3) \mapsto (t^4g_2,t^6g_3).$$
The coarse moduli space is
$$\oM_{1,1}=\oA_1 \simeq \Proj k[g_2,g_3]\simeq \bP^1.$$
The symmetric square of the stack has the following
quotient-stack presentation:
\begin{eqnarray*}
(\ocM_{1,1} \times \ocM_{1,1})/\fS_2
&=&\left(\Spec k[g_2,g_3,h_2,h_3]-Z\}\right/H \\
Z&=&\{(g_2,g_3,h_2,h_3): g_2=g_3=0 \text{ or } h_2=h_3=0\}.
\end{eqnarray*}
Here $H$ is the group generated by the torus
$$(t,u)\cdot(g_2,g_3,h_2,h_3)\mapsto 
(t^4g_2,t^6g_3,u^4h_2,u^6h_3)$$
and the involution
$$(g_2,g_3,h_2,h_3) \mapsto (h_2,h_3,g_2,g_3),$$
i.e., $H=\fS_2 \ltimes \bG_m^2$ where $\fS_2$
acts on $\bG_m^2$ by interchanging the factors.  

The coarse moduli space of the stack is the 
invariant-theory quotient for the action of $H$.  
Consider the elements
$p\in k[g_2,g_3,h_2,h_3]$ with the following properties:
\begin{enumerate}
\item{$p(g_2,g_3,h_2,h_3)=p(h_2,h_3,g_2,g_3)$;}
\item{$p(t^2g_2,t^3g_3,u^2h_2,u^3h_3)=(tu)^N p(g_2,g_3,h_2,h_3)$
for some $N$.}
\end{enumerate}
This ring is generated by $g_2h_2,g_3h_3,$ and $g_2^3h_3+g_3^2h_2^3$
and 
$$\Proj k[g_2h_2,g_3h_3,g_2^3h_3^2+g_3^2h_2^3] \simeq \bP(4,6,12)
\simeq \bP(2,3,6),$$
which explains why the weights of $b$ are $(2,3,6)$.
See \cite{Ig2}, Theorem 3, for a discussion
in terms of the modular forms for 
$\mathrm{Sp}(2,\bZ)$ (see Remark~\ref{rem:Igusa}).  

\noindent {\em proof of proposition:}
The first equation follows because $b:\widetilde{X}\ra X$
has weights $(2,3,6)$.  
As for the second,
$\widetilde{\Xi}$ is the proper transform of the divisor
$\{G=0\} \subset X$ parametrizing forms
admitting an bielliptic involution.  When a bielliptic curve
of genus two specializes to a stable curve in $\Delta_1$,
the bielliptic involution specializes to a morphism
exchanging the elliptic components.  Therefore,
$\widetilde{\Xi} \cap \widetilde{\Delta}_1 \subset 
\widetilde{\Delta}_1 $
is the diagonal in the symmetric square, which is
cut out by a form of weighted-degree twelve. 

For the third equation, $\widetilde{\Delta}_0$ is
the proper transform of $\{D=0\}$.  The intersection
$\widetilde{\Delta}_0 \cap \widetilde{\Delta}_1 \subset 
\widetilde{\Delta}_1$
is the locus where the discriminant
$\Delta=g_2^3-27g^2_3$
vanishes, and thus has weighted-degree six.  
The last equation follows because $G$ has weighted
degree thirty and $D$ has weighted degree ten.$\square$

One consequence of this analysis is worth mentioning.  
\begin{prop}\label{prop:blowpol}
For sufficiently small $\epsilon>0$, the divisor
$$\widetilde{\Delta}_0+(6-\epsilon)\widetilde{\Delta}_1$$
is ample on $\widetilde{X}$.  
\end{prop}
{\em proof:}
The divisor can be expressed
$$
b^*(\text{ample divisor})-\epsilon(\text{$b$-exceptional divisor}),$$
which yields a polarization of the blow-up $b:\widetilde{X}\ra X$.
$\square$

\section{Birational geometry}
\subsection{Divisor classes  and birational contractions of
$\oM_2$}
\label{subsect:canM2}
It is well known that the rational divisor class group of $\oM_2$
is freely generated by the boundary divisors
$\Delta_0$ and $\Delta_1$;  
see \cite{HaMo} for a nice account of divisors
on $\oM_g$ for arbitrary $g$.  When $g=2$,
Proposition~\ref{prop:bircont} gives an elementary proof of
this fact:  Since $X\simeq \bP(1,2,3,5)$ its divisor class
group has rank one and is generated by the discriminant 
divisor $\{D=0\}$; the same holds
true for $X-p$.  By Proposition~\ref{prop:codimone}
$\oM_2$ is isomorphic
to $\widetilde{X}$ up to codimension $\le 1$, so these have 
isomorphic class groups.  It follows that the rational
divisor class group of $\oM_2$ is generated by $\Delta_1$
and the proper transform of the discriminant, which is 
just $\Delta_0$.  

The nef cone of $\oM_2$ is also well-known.  
We will not give a self-contained proof here, but rather rely on the
general result of Cornalba-Harris \cite{CH}:
\begin{thm} \label{thm:CH}
The line bundle $a\lambda -b\delta$ is nef on $\oM_g,g\ge 2,$ if and only if
$a\ge 11b \ge 0$.  
\end{thm}
Here $\lambda$ is the pull-back of the polarization on $\oA_g$ via
the extended
Torelli map $\overline{t}:\oM_g \ra \oA_g$ (see \S\ref{subsect:comp}).

To apply this in our situation, we observe that
$$\lambda\equiv \frac{1}{10}(\Delta_0+\Delta_1)$$
over $\oM_2$ (see \cite{HaMo} pp. 175).
The factor $10$ can be
explained by the fact that $\overline{t}(\Delta_1)$ is
defined by the vanishing of a modular form of weight ten
(see Remark~\ref{rem:Igusa} and \cite{Ig2}).  Substitution
gives the first part of
\begin{prop}\label{prop:nef2}
The nef cone of $\oM_2$ is generated by the divisors
$\Delta_0+\Delta_1$ and $\Delta_0+6\Delta_1$,
respectively.  
These are both semiample, inducing 
the birational contractions
$$\overline{t}:\oM_2 \ra \oA_2 \quad \beta:\oM_2 \ra X.$$
The rational map $\gamma:\oM_2 \dashrightarrow \widetilde{X}$
is an isomorphism.
\end{prop}
{\em remainder of proof:}
Of course, $\Delta_0+\Delta_1$ is semiample and induces the
birational contraction morphism
$\overline{t}:\oM_2 \ra \oA_2$.  As
$\widetilde{X}$ and $\oM_2$ are isomorphic in codimension one,
Proposition~\ref{prop:blowpol} says that 
$\widetilde{\Delta}_0+(6-\epsilon)\widetilde{\Delta}_1$
is ample on $\widetilde{X}$ and the corresponding
divisor $\Delta_0+(6-\epsilon)\Delta_1$ is ample on $\oM_2$.
It follows that $\widetilde{X}$ and $\oM_2$ are each
isomorphic to $\Proj$ of 
$$\oplus_{n\ge 0} \Gamma(\cO_{\widetilde{X}}(
n(\widetilde{\Delta}_0+(6-\epsilon)\widetilde{\Delta}_1)))
\simeq  \oplus_{n\ge 0} \Gamma(\cO_{\oM_2}(
n(\Delta_0+(6-\epsilon)\Delta_1))) 
.$$
In particular, the rational map $\gamma$ is an isomorphism.
Thus the contractions $b:\widetilde{X}\ra X$ and 
$\beta:\oM_2 \dashrightarrow X$ coincide.  $\square$

\subsection{Canonical class of $\oM_2$}

The canonical class $K_{\oM_2}$ can also be computed by
elementary methods.  We know that
$$\omega_{\bP(1,2,3,5)}\simeq \cO_{\bP(1,2,3,5)}(-11),$$
which also follows from the third
part of Proposition~\ref{prop:ramq}.  Since the discriminant
has degree ten, we find
$$K_X=-\frac{11}{5}\{D=0\}.$$
Applying the formulas of Proposition~\ref{prop:symsquare}, we obtain
$$K_{\widetilde{X}}\equiv
-\frac{11}{5}\widetilde{\Delta}_0-\frac{16}{5}\widetilde{\Delta}_1$$
$$K_{\widetilde{X}}+\frac{1}{2}\widetilde{\Delta}_1+\frac{1}{2}\widetilde{\Xi}
\equiv
-\frac{7}{10}\widetilde{\Delta}_0+\frac{3}{10}\widetilde{\Delta}_1.$$
Since $\widetilde{X}$ and $\oM_2$ agree in codimension one, they have
the same canonical class  
$$K_{\oM_2}\equiv -\frac{11}{5}\Delta_0-\frac{16}{5}\Delta_1.$$
In particular, we obtain
\begin{equation} \label{logcancomp}
K_{\oM_2}+ \alpha \Delta_0+
\frac{1+\alpha}{2}\Delta_1+\frac{1}{2}\Xi\equiv 
(-\frac{7}{10}+\alpha)\Delta_0
+(\frac{3}{10}+\alpha/2)\Delta_1.
\end{equation}
\begin{rem}
The importance of this divisor stems from the ramification
equation of Proposition~\ref{prop:ramify}.  This divisor class
pulls back to the class
$$K_{\ocM_2}+\alpha\delta$$
on the moduli stack.  We shall interpret log canonical
models of the moduli stack using this divisor.  
\end{rem}

Using the computations of Proposition~\ref{prop:symsquare},
we obtain a discrepancy equations for $\beta:\oM_2\ra X$
\begin{equation}
\label{M2discrep1}
K_{\oM_2}+\alpha\Delta_0+\frac{1}{2}\Xi=
\beta^*(K_X+\alpha\{D=0\}+\frac{1}{2}\{G=0\})+(4-6\alpha)\Delta_1
\end{equation}
\begin{equation}
\label{M2discrep2}
K_{\oM_2}+\alpha\Delta_0+\frac{1+\alpha}{2}\Delta_1+\frac{1}{2}\Xi=
\beta^*(K_X+\alpha\{D=0\}+\frac{1}{2}\{G=0\})+\frac{9-11\alpha}{2}\Delta_1.
\end{equation}

\subsection{Generalities on log canonical models}
See \cite{FA} \S 2
for definitions of the relevant terms and technical background.  
Let $V$ be a normal projective variety, $D=\sum_{i=1}^n a_iD_i$
a $\bQ$-divisor such that $0\le a_i\le 1$ and $K_V+D$
is $\bQ$-Cartier.  Abusing notation, we write $V-D$ 
for $V-\cup_i D_i$.  
\begin{defn} \label{defn:strict}
$(V,D)$ is a {\em strict log canonical
model} if 
$K_V+D$ is ample, $(V,D)$ has log canonical singularities, and 
$V-D$ has canonical singularities.
\end{defn}
The idea here is to realize $(V,D)$ as a log
canonical model without introducing boundary divisors over
$V-D$. This is natural if we want to 
respect the geometry of the open complement.  

The following recognition criterion for
strict log canonical models is based on 
\cite{FA} \S 2:
\begin{prop}\label{prop:strict}
Consider birational projective contractions
$\rho:\widetilde{V}\ra V$ where the exceptional
locus of $\rho$ is divisorial, 
$D'_i$ denotes the proper transform of $D_i$, and $E_j$ (resp. $F_k$)
denotes the exceptional divisors of $\rho$ with $\rho(E_j)\subset D$ (resp.
$\rho(F_k) \not \subset D$).

The following are equivalent:
\begin{enumerate}
\item{$(V,D)$ is a strict log canonical model.}
\item{For some
resolution of singularities $\rho:\widetilde{V} \ra V$, with
the union of the exceptional locus and 
$\cup_i D'_i$ normal crossings, there exist
$b_j\in \bQ\cap [0,1]$ so that
$$\widetilde{D}=\sum_i a_i D'_i + \sum_j b_j E_j$$
satisfies the formula
\begin{equation}
\label{eq:discrep}
K_{\widetilde{V}}+\widetilde{D} \equiv \rho^*(K_V + D)+
\sum_j d_j E_j + \sum_k e_k F_k, \quad d_j,e_k \ge 0.
\end{equation}
For each such choice of $b_j$ 
$$K_{\widetilde{V}}+\widetilde{D} - \sum_j d_j E_j - \sum_k e_k F_k$$
is semiample and induces $\rho$.
Furthermore, we may take the 
$b_j=1$.}
\item{For some contraction
$\rho:\widetilde{V} \ra V$,
there exist $b_j\in \bQ\cap [0,1]$ so that
$$\widetilde{D}=\sum_i a_i D'_i + \sum_j b_j E_j$$
satisfies the formula
$$
K_{\widetilde{V}}+\widetilde{D} \equiv \rho^*(K_V + D)+
\sum_j d_j E_j + \sum_k e_k F_k, \quad d_j,e_k \ge 0,
$$
$(\widetilde{V},\widetilde{D})$ is log canonical, and 
$\widetilde{V}-\widetilde{D}$ has canonical singularities.
The divisor
$$K_{\widetilde{V}}+\widetilde{D} - \sum_j d_j E_j - \sum_k e_k F_k$$
is semiample and induces $\rho$.}
\end{enumerate}
\end{prop}
Some general facts are worth mentioning before we indicate
the proof.  First, we can decide
whether a pair is canonical or log canonical by computing
discrepancies on {\em any} resolution.  Second, discrepancies
increase as the coefficients of the log divisor are
decreased \cite{FA} 2.17.3.  Third, in situations
(2) and (3) the pair $(V,D)$ is the {\em log canonical model} of 
$(\widetilde{V},\widetilde{D})$;  such
models are unique \cite{FA} 2.22.1.  

\noindent
{\em proof:} 
It is trivial that the second statement implies the third.  
To see that the third implies the first, take a resolution
for $(\widetilde{V},\widetilde{D})$ so that the union of the
exceptional locus and all the proper transforms of the
$D'_i,E_j,$
and $F_k$ is normal crossings.  
Comparing discrepancies for 
$(\widetilde{V},\widetilde{D})$ and $(V,D)$, using the fact
the coefficients of components in $\widetilde{D}$ are
at least as large as the coefficents of the corresponding
components appearing in $\rho^*(K_V+D)$, we find 
that $(V,D)$ is log-canonical 
and has canonical singularities 
along $V-D$.
Since $\rho^*(K_V+D)$ induces $\rho$,
$K_V+D$ must be ample on $V$.  

For the remaining implication, since $(V,D)$ has log canonical
singularities and canonical singularities away from $D$,
the discrepancy equation (\ref{eq:discrep}) follows.  Since
$K_V+D$ is ample on $V$, its pull-back to $\widetilde{V}$ is
semiample and induces $\rho$.  
$\square$

We shall also need the following basic fact,
a special case of \cite{FA} 20.2, 20.3:
\begin{prop}\label{prop:ramifylc}
Let $W$ be a smooth variety 
and $h:W\ra V$ be a finite dominant morphism to a normal variety.
Let $D=\sum_i a_i D_i, 0\le a_i \le 1$ be a $\bQ$-divisor on $V$ 
containing all the divisorial components of the branch locus of $h$.
Let $\bar{D}$
be a $\bQ$-divisor
on $W$ so that 
$\supp (h^{-1}(D))=\supp (\bar{D})$ and
$h^*(K_V+D)=K_W+\bar{D}$.
Then $(V,D)$ has log canonical singularities along $D$ iff
$(W,\bar{D})$ has log canonical singularities along $\bar{D}$.  

If $D$ and $\bar{D}$ have multiplicity one at each component then
$h^*(K_V+D)=K_W+\bar{D}$ follows from the other assumptions.
\end{prop}

\subsection{Example of $\oM_g$}
The standpoint of this section owes a great deal
to Mumford \cite{Mu1} \cite{Mu2}:
\begin{thm} \label{thm:Mumford}
For $g\ge 4$, the pair
$(\oM_g,\Delta)$
is a strict log canonical model.  
\end{thm}
\begin{rem}\label{rem:stacklc}
This is also
the natural log canonical model
from the point of view of the moduli stack.
Indeed, for $g\ge 4$
the locus in $M_g$ of curves with automorphisms has codimension
$\ge 2$, so the branch divisor of $Q:\ocM_g \ra \oM_g$ is just
$\Delta_1$;  over $\Delta_1$, we have simple ramification.
We therefore have (\cite{HaMu} pp. 52)
$Q^*K_{\oM_g}=K_{\ocM_g}-\delta_1$ and thus
$$Q^*(K_{\oM_g}+\Delta)=K_{\ocM_g}+\delta.$$
\end{rem}
{\em sketch proof:}
We first check that $K_{\oM_g}+\Delta$ is ample.  The formula
from \cite{HaMu} \S 2 (or \cite{HaMo})
$$K_{\oM_g}=13\lambda-2\Delta_0-3/2\Delta_{1}-2\sum_{2\le i\le g/2}\Delta_i$$
gives
$K_{\oM_g}+\delta =13\lambda-\delta$, which is ample by Theorem~\ref{thm:CH}
(see also \cite{Mu1}). 

The singularity analysis follows \cite{HaMu}.
$\oM_g$ has canonical singularities by Theorem 1 of \cite{HaMu}.
To show that $(\oM_g,\Delta)$ has log canonical singularities,
we use the fact that $\oM_g$ is \'etale-locally
a quotient of a smooth
variety by a finite group.   At $[C]$ it has a local presentation
$$h:T_{[C]}\ocM_g \twoheadrightarrow T_{[C]}\ocM_g/\Aut(C)$$
in terms of its tangent space
$$T_{[C]}\ocM_g=\mathrm{Ext}^1(\Omega^1_C,\cO_C).$$

We analyze the quotient morphism using 
Proposition~\ref{prop:ramifylc}.
The preimage of the boundary divisor corresponds to a union of hyperplanes 
$\widehat{\delta}\subset T_{[C]}\ocM_g$, meeting in normal crossings.
The pair $(T_{[C]}\ocM_g,\bar{\Delta})$ then has log canonical singularities.  
An application of Proposition~\ref{prop:ramifylc}, utilizing 
the ramification discussion in Remark~\ref{rem:stacklc}, implies
$(K_{\ocM_g},\Delta)$ is log canonical.  $\square$

\begin{rem} \label{rem:fullforce}
Using the full force of Theorem~\ref{thm:CH} we get a sharper statement.
Consider the pair
$$(\oM_g,\alpha\Delta_0 +\frac{1+\alpha}{2}\Delta_1 + 
\alpha(\Delta_2+\ldots+\Delta_{\lfloor g/2 \rfloor})),$$
with log canonical divisor pulling back to
$K_{\ocM_g}+\alpha \delta$ on the moduli stack.  
The $\bQ$-divisor $K_{\ocM_g}+\alpha\delta$ is the pull-back
of an ample line bundle if and only if $9/11<\alpha\le 1$.  
Since $\oM_g$ is a locally a quotient of a smooth variety
by a finite group, all divisors on $\oM_g$ are $\bQ$-Cartier.
An easy computation with the discrepancy equation (\ref{eq:discrep})
then shows that
the pair remains log canonical even as the coefficients
are reduced.
\end{rem}

\subsection{Application to $\ocM_2$}
The source of the special difficulties
in this case is the fact that
$K_{\oM_2}+\Delta$ is not effective.  
Indeed, in \S\ref{subsect:canM2} we computed
$$K_{\oM_2}\equiv -\frac{11\Delta_0+6\Delta_1}{5}$$
so $K_{\oM_2}+\Delta=-\text{effective divisor}$.  

In order to recover a result analogous to Theorem~\ref{thm:Mumford},
we must take the `log canonical model of the moduli stack',
as interpreted on $\oM_2$ via Proposition~\ref{prop:ramify}:
\begin{thm} \label{thm:main}
Consider the log canonical model of $\ocM_2$ with respect to the 
$K_{\ocM_2}+\alpha\delta$, i.e., 
the log canonical model of $\oM_2$ with
respect to 
$$K_{\oM_2}+\alpha \Delta_0+\frac{1+\alpha}{2}\Delta_1+\frac{1}{2}\Xi.$$
\begin{enumerate}
\item{For $9/11<\alpha\le 1$, we recover $\oM_2$.}
\item{For $7/10 < \alpha \le 9/11$ we recover the invariant theory
quotient $X\simeq \bP(1,2,3,5)$.}
\item{For $\alpha =7/10$ we get a point; the log canonical divisor
fails to be effective for $\alpha<7/10$.}
\end{enumerate}
\end{thm}
{\em proof:}
The necessary ampleness results have already been stated.
Proposition~\ref{prop:nef2} and Equation (\ref{logcancomp}) 
imply the log canonical divisor on $\oM_2$ is ample if and only if
$\alpha>9/11$.  Proposition~\ref{prop:ramq} implies that
$$K_X+1/2\{G=0\}+\alpha\{D=0\}$$
is positive on $X$ if and only if $\alpha>7/10$.  
When $\alpha=7/10$ it is
zero and when $\alpha<7/10$ it is negative.  

It remains to verify the singularity conditions:
First, we check that $\oM_2$ has canonical singularities away
from $\Delta_0,\Delta_1,$ and $\Xi$.  
Suppose that $C$ is not in the boundary and
does not admit admit a bielliptic involution.
In Proposition~\ref{prop:m2} we saw
$$M_2 \simeq X-\{D=0\} \subset \bP(1,2,3,5),$$
so we need to analyze the singularities of $\bP(1,2,3,5)-\{D=0\}$.  
A point in weighted projective space is nonsingular when the weights
corresponding to its non-vanishing coordinates are relatively prime,
so the only possible singularity occurs when $A=B=C=0$.  
The corresponding binary sextic form is
$$x(x^5+y^5),$$
the unique form with an automorphism group of order five 
\cite{Bo} pp. 51 \cite{Ig1} pp. 645.  At this point, $\bP(1,2,3,5)$
is locally isomorphic to the cyclic quotient singularity
$\frac{1}{5}(1,2,3)$, i.e., the quotient of $\bA^3$ under the action
$$(a,b,c) \mapsto (\zeta a, \zeta^2 b, \zeta^3 c) \quad 
\zeta\neq 1 \in \mu_5.$$
This is canonical by the Reid-Tai criterion;  see \cite{HaMu}
pp. 28 for a general result.

Second, we address the singularities along the boundary.
For $\alpha>9/11$ we need that
$$K_{\oM_2}+\alpha \Delta_0 + \frac{1+\alpha}{2}\Delta_1+1/2\Xi$$
is log canonical.  When $\alpha\le 9/11$,
Proposition~\ref{prop:strict}
and the discrepancy computation (\ref{M2discrep2}) reduce us to showing 
that this is log canonical.    
Since $\oM_2$ is $\bQ$-factorial and discrepancies
increase as coefficients of log divisors decrease \cite{FA} 2.17.3,
it suffices to verify that  
\begin{equation}
\label{logmain}
K_{\oM_2}+\Delta_0 + \Delta_1+\frac{1}{2}\Xi
\end{equation}
is log canonical.  

The proof relies on the description of the boundary divisors
in terms of the local presentation 
$$T_{[C]}\ocM_2/\Aut(C),$$
as sketched in \S\ref{subsect:stackgeo}.   The key observation is that
$\Xi$ does not play a r\^ole in the analysis.
Each bielliptic involution
acts on $T_{[C]}\ocM_2$ by reflection across the corresponding hyperplane 
in $\widehat{\xi}$, so the quotient 
$$h:T_{[C]}\ocM_2 \ra T_{[C]}\ocM_2/H, 
\quad H=\Aut(C)/\left<\iota\right>,$$
has simple ramification along $\widehat{\xi}$.  Since $\Xi$ has
coefficient $1/2$ in (\ref{logmain}), $\widehat{\xi}$ 
does not appear in the pull-back of the
log canonical divisor to $T_{[C]}\ocM_2$.  

Thus (\ref{logmain}) pulls back to 
$$K_{T_{[C]}\ocM_2}+\widehat{\delta},$$ 
and we have seen that $\widehat{\delta}$
is normal crossings.  It follows that $(T_{[C]}\ocM_2,\widehat{\delta})$
is log canonical and
Proposition~\ref{prop:ramifylc} gives the desired result.$\square$

\

\noindent
Rice University, MS 136\\
Department of Mathematics \\
PO Box 1892\\
Houston, TX 77251-1892\\
USA \\
\texttt{hassett@math.rice.edu}


\begin{thebibliography}{}
\bibitem[ACV]{ACV}
D. Abramovich, A. Corti, and  A. Vistoli,
Twisted bundles and admissible covers,
{\em Comm. Algebra} {\bf  31} (2003), no. 8, 3547--3618.



\bibitem[Bo]{Bo}
O. Bolza, On binary sextics with linear transformations
into themselves, {\em Amer. J. Math.} {\bf 10} (1887)
47-70.  


\bibitem[Ca]{Ca}
A. Cayley, Tables for the binary sextic, {\em Amer. J.
 Math.} {\bf 4} (1881), 379--384.  Reprinted in:
{\em Collected Mathematical Papers}, vol. XI, pp. 372-376,
Cambridge University Press, 1896.  

\bibitem[Cl]{Cl}
A. Clebsch, {\em Theorie der bin\"aren algebraischen Formen},
Verlag von B.G. Teubner, Leipzig, 1872.  

\bibitem[CH]{CH}
M. Cornalba and J. Harris, 
Divisor classes associated to families of stable varieties, 
with applications to the moduli space of curves,
{\em Ann. Sci. \'Ecole Norm. Sup. (4) }{\bf 21 }
(1988), no. 3, 455--475.
 


\bibitem[El]{El}
E.B. Elliott, {\em An Introduction to the Algebra of Quantics},
Oxford University/Clarendon Press, 1895.

\bibitem[FA]{FA}
J. Koll\'ar, {\em Flips and abundance
for algebraic threefolds}, {\em Ast\'erisque}
{\bf 211}, 1992.

\bibitem[HaMo]{HaMo}
J. Harris and I. Morrison, {\em Moduli of Curves}, 
Springer-Verlag, New York, 1998.  

\bibitem[HaMu]{HaMu}
J. Harris and D. Mumford,  On the Kodaira dimension
of the moduli space of curves, {\em Invent. math.}
{\bf 67} (1982), 23--86.  

\bibitem[Ha]{Ha}
R. Hartshorne,
{\em Algebraic geometry},
Springer-Verlag, New York-Heidelberg, 1977. 

\bibitem[Hi]{Hi}
D. Hilbert, {\em Theory of algebraic invariants},
Cambridge University Press, 1993.  

\bibitem[Ig1]{Ig1}
J. Igusa, 
Arithmetic variety of moduli for genus two,
{\em Ann. of Math. (2)} {\bf 72} (1960) 612--649.

\bibitem[Ig2]{Ig2}
J. Igusa, 
 On Siegel modular forms of genus two,
{\em Amer. J. Math.}  {\bf 84} (1962) 175--200.

\bibitem[Ig3]{Ig3}
J. Igusa, 
A desingularization problem in the theory of Siegel modular functions,
{\em Math. Ann.}  {\bf 168} (1967) 228--260

\bibitem[LM]{LM}
G. Laumon and L. Moret-Bailly,
{\em Champs alg\'ebriques}, 
Springer-Verlag, Berlin-Heidelberg, 2000.

\bibitem[Mu1]{Mu1}
D. Mumford, 
Stability of projective varieties,
{\em Enseignement Math. (2)} {\bf 23} (1977), no. 1-2, 39--110.

\bibitem[Mu2]{Mu2}
D. Mumford, 
Hirzebruch's proportionality theorem in the noncompact case,
{\em Invent. Math.}{\bf 42} (1977), 239--272.

\bibitem[GIT]{GIT}
D. Mumford, J. Fogarty, and F. Kirwan,
{\em Geometric Invariant Theory}, third enlarged edition,
Springer-Verlag, Berlin-Heidelberg, 1994.

\bibitem[Na1]{Na1}
Y. Namikawa, 
On the canonical holomorphic map from the 
moduli space of stable curves to the Igusa monoidal transform,
{\em Nagoya Math. J.} {\bf 52} (1973), 197--259.

\bibitem[Na2]{Na2}
Y. Namikawa, 
A new compactification of the Siegel space and degeneration of 
Abelian varieties I/II 
{\em Math. Ann.} {\bf 221} (1976), no. 2, 97--141/ no. 3, 201--241.




\bibitem[Sa]{Sa} G. Salmon, {\em Lessons introductory to the 
modern higher algebra,} 
2nd edition, 
Hodges Smith, Dublin, 1866.  

\bibitem[SB]{SB} N. I. Shepherd-Barron, Canonical rings for moduli spaces
of abelian varieties, preprint (2003).

\bibitem[Sc]{Sc} I. Schur, {\em Vorlesungen \"uber Invariententheorie},
Bearbeitet und herausgegeben von Helmut Grunsky,
Springer-Verlag, Berlin-New York 1968

\bibitem[Sp]{Sp}
T.A. Springer, 
{\em Invariant theory},
Lecture Notes in Mathematics, Vol. 585,
Springer-Verlag, Berlin-New York, 1977.

\bibitem[Vi]{Vi}
A. Vistoli, 
The Chow ring of $\cM_2$, Appendix to 
"Equivariant intersection theory" by  D. Edidin and W. Graham,
Invent. Math. {\bf 131} (1998), no. 3, 635--644.

\end{thebibliography}
\end{document}